\documentclass[11pt]{amsart}
\usepackage{amsmath,amssymb, verbatim, amsthm}

\numberwithin{equation}{section}

\newtheorem{thm}{Theorem}[section]
\newtheorem{prop}[thm]{Proposition}
\newtheorem{lemma}[thm]{Lemma}
\newtheorem{cor}[thm]{Corollary}
\newtheorem{example}[thm]{Example}
\newtheorem{remark}[thm]{Remark}
\newtheorem{definition}[thm]{Definition}

\newenvironment{ex}{\begin{example}\rm}{\end{example}}
\newenvironment{rem}{\begin{remark}\rm}{\end{remark}}

\DeclareMathOperator{\coop}{coop}                                              
\DeclareMathOperator{\op}{op}

\newcommand{\DOT}{\setlength{\unitlength}{1pt}\begin{picture}(2.5,2)
               (1,1)\put(2,3.5){\circle*{3}}\end{picture}}
\newcommand{\db}{\DOT_{\beta}}
\renewcommand{\_}[1]{\mbox{$_{\left( #1 \right)}$}}   
\renewcommand{\k}{{\mathbb K}}
\newcommand{\com}{\Delta}
\newcommand{\Z}{{\mathbb Z}}
\newcommand{\id}{\mbox{\rm id\,}}      
\newcommand{\Hom}{\mbox{\rm Hom\,}}
\newcommand{\w} {\omega}

\newcommand{\ursn}{{\mathfrak{u}}_{r,s}({\mathfrak{sl}}_n)}

\newcommand{\baru}{\overline{\mathfrak{u}_{r,s}({\mathfrak{sl}}_n)}}
\newcommand{\E}{\mathcal E}
\newcommand{\F}{\mathcal F}
\newcommand{\D}{\mathcal D}
\renewcommand{\ker}{\mbox{\rm Ker\,}}  

\newcommand{\Id}{\mathop{\rm Id}}
\renewcommand{\Im}{\mathop{\rm Im\,}}

\title[Factorization of simple modules]{
Factorization of simple modules for certain \\
  pointed Hopf Algebras}
\author{Mariana Pereira}
\address{Centro de Matem\'atica\\
        Facultad de Ciencias,
        Igu\'a 4225 \\
        Montevideo, 11400 \\
        Uruguay}
\email{mariana@cmat.edu.uy}
\date{21 Feb 2007}

\begin{document}
\begin{abstract}
We study the representations of two types of pointed Hopf algebras: restricted two-parameter quantum groups, and the Drinfel'd doubles of rank one pointed Hopf algebras of nilpotent type. We study, in particular, under what conditions a simple module can be factored as the tensor product of a one dimensional module with a module that is naturally a module for the quotient by central group-like elements. For restricted two-parameter quantum groups, given $\theta$ a primitive $\ell$th root of unity, the factorization of simple $\mathfrak{u}_{\theta^y, \theta^z}(\mathfrak{sl}_n)$-modules is possible, if and only if $\gcd((y-z)n, \ell)=1$. For rank one pointed Hopf algebras, given the data $\D=(G,\chi, a)$, the factorization of simple $D(H_{\D})$-modules is possible if and only if $|\chi(a)|$ is odd and $|\chi(a)|=|a|=|\chi|$. Under this condition, the tensor product of two simple $D(H_{\D})$-modules is completely reducible, if and only if the sum of their dimensions is less than or equal to $|\chi(a)|+1$.

\end{abstract}

\maketitle

\section*{Introduction} 
In 2003, Radford introduced a new method to construct simple left modules for the Drinfel'd double of a graded Hopf algebra with commutative and cocomutative bottom term \cite{radford}. Until then, simple modules for such algebras were usually constructed by taking quotients of Verma modules by maximal submodules. This new method gives a more explicit construction, in the sense that the simple modules are given as subspaces of the Hopf algebra itself. With this construction, simple modules for the Drinfel'd double of a Hopf algebra are parametrized by group-like elements of the Drinfel'd double. 
Recently, new examples of pointed Hopf algebras have been found (e.g \cite{krop-radford, a-s-1}) and we would like to understand their modules. Our approach is, if possible, to reduce the problem to that of understanding modules for special types of Hopf algebras, particularly when these special types are well known (e.g traditional quantum groups, Taft algebras). The reduction is by factorizing simple modules as the tensor product of a one-dimensional module with a module that is naturally a module for the quotient of the algebra by the central group-like elements. We study this factorization, using Radford's method, for two types of pointed Hopf algebras: restricted two-parameter quantum groups, and the Drinfel'd doubles of rank one pointed Hopf algebras of nilpotent type; and we find necessary and sufficient conditions for the factorization to be possible.

In Section \ref{prelims} we give the basic definitions and present Radford's method to construct simple modules for $D(H)$ when $H$ is a finite-dimensional graded Hopf algebra with a commutative group algebra as its bottom term. These modules are in one to one correspondence with $G(H^*)\times G(H)$, the product of the set of group-like elements of the dual of $H$ with the set of group-like elements of $H$. In Proposition \ref{tensor products} we show that given $(\beta, g)$ and $(\beta', g')$ $\in G(H^*)\times G(H)$, if the tensor product of the corresponding simple $D(H)$-modules is also simple, then it is isomorphic to the module corresponding to $(\beta \beta', g g')$.

In Section \ref{quantum}, we recall the definition of the restricted two-parameter quantum groups $\ursn$. Under certain conditions on the parameters $r$ and $s$, $\ursn$ is a Drinfel'd double and Radford's method can be used to construct its simple modules. In Proposition \ref{1dim-prop} we characterize the group-like elements of $\ursn$ that give one-dimensional $\ursn$-modules. In Proposition \ref{centralu} we describe the central group-like elements of $\ursn$ and in Theorem \ref{final-ursn} we give necessary and sufficient conditions on the parameters $r$ and $s$ for any simple $\ursn$-module to factor as the tensor product of a one-dimensional module with a $\ursn$-module that is also a module for a quotient of $\ursn$.

In Section \ref{pointed}, we present the rank one pointed Hopf algebras of nilpotent type, $H_{\D}$ defined by Krop and Radford \cite{krop-radford}, and get analogous results for simple modules for their Drinfel'd doubles as those for $\ursn$ in Proposition \ref{1dimK-R} and Theorem \ref{decomposition-KR}. In Theorem \ref{last-pointed} we use these results to study the reducibility of the tensor product of two simple $D(H_{\D})$-modules. 
\section{Preliminaries}\label{prelims}
In what follows $\k$ is an algebraically closed field of characteristic 0. All vector spaces and tensor products are over $\k$. For any vector space $V$, if $f \in V^*$ and $v\in V$, the evaluation of $f$ on $v$ will be denoted by $\langle f, v\rangle$. If $(C, \Delta, \epsilon)$ is a coalgebra, usually denoted by $C$, and $c \in C$  we write $ \com (c) = \sum c\_1 \otimes c\_2 $ (Heyneman-Sweedler notation). The opposite coalgebra $(C, \Delta^{\coop}, \epsilon)$, usually denoted by $C^{\coop}$, has comultiplication given in terms of the comultiplication of $C$: $\Delta^{\coop}(c) = \sum c\_2 \otimes c\_1$, for all $c \in C$. If $(M, \delta)$ is a right $C$-comodule and $m \in M$, we write $\delta(m)= \sum m\_0 \otimes m\_1 \in M \otimes C$. We will say comodule for right comodule and module for left module. The coalgebra $C$ is a $C^*$-module with action $f \cdot c = \sum c\_1\langle f, c\_2 \rangle$ for all $f \in C^*$ and $c \in C$. 

The set of group-like elements of $C$ is $G(C)=\{c \in C: \Delta(c) = c \otimes c \mbox{ and } \epsilon(c)=1 \}$, and it is a linearly independent set. If $H$ is Hopf algebra, then $G(H)$ is a group and its linear span, $\k G(H)$, is a Hopf algebra. If $A$ is a finite-dimensional algebra, then $G(A^*)=\mbox{Alg}(A, \k)$, the set of algebra maps from $A$ to $\k$. If $C$ is a coalgebra, then $\mbox{Alg}(C,\k)$ is an algebra with the {\em convolution product}: $\langle fg , c \rangle = \sum \langle f, c\_1 \rangle \langle g , c\_2 \rangle$ for all $f,g \in C^*$ and $c \in C$.

Let $H$ be a finite-dimensional Hopf algebra over $\k$
with antipode $S$. The {\em Drinfel'd double of $H$}, $D(H)$, is
$$D(H)=(H^*)^{\coop}\otimes H$$ as a coalgebra. The algebra structure is given by
$$ (g\otimes h)(f\otimes k) = \sum g \left( h\_1 \rightharpoonup f \leftharpoonup S^{-1}(h\_3) \right) \otimes h\_2k,$$
for all $g,f \in H^*$ and $h,k \in H$; where $\langle a \rightharpoonup f , b\rangle = \langle f, ba \rangle $ and $\langle f \leftharpoonup a , b \rangle =\langle f ,ab \rangle$, for all $a, b \in H$ and $f \in H^*$. With this structure and antipode $\mathcal{S}$ given by $\mathcal{S}(f \otimes h) = (\epsilon \otimes S(h))((f \circ S^{-1})\otimes 1)$, $D(H)$ is a Hopf algebra. 

For any bialgebra $H$, a {\em left-right Yetter-Drinfel'd module} is a triple $(M, \cdot, \delta)$ where $(M, \cdot)$ is a left $H$-module and $(M, \delta)$ is a right $H$-comodule, and the following compatibility condition is satisfied:
$$\sum h\_1 \cdot m\_0 \otimes h\_2 m\_1 = \sum ( h\_2\cdot m)\_0 \otimes (h\_2 \cdot m)\_1 h\_1.$$
The category of left-right Yetter-Drinfel'd modules over a bialgebra $H$ will be denoted by ${}_H {\mathcal YD}^H$.

\begin{prop}[Majid \cite{majid}]\label{yetter} Let H be a finite-dimensional Hopf algebra. Then
$D(H)$-modules are left-right Yetter-Drinfel'd modules and conversely. 
Explicitly, if $M$ is a left-right Yetter-Drinfel'd module, then it is a $D(H)$-module with the same action of $H$ and the action of $H^*$ given by
\begin{equation}\label{action-dual-general}
f \cdot m = \sum \langle f,m\_1 \rangle m\_0,
\end{equation}
for all $f$ in $H^*$ and $m$ in $M$.
\end{prop}

If $M,N \in {}_H {\mathcal YD}^H$, $M \otimes N$ is also Yetter-Drinfel'd module over $H$. The Yetter-Drinfel'd structure is given by the action
$$h \cdot (m \otimes n) =\sum h\_1 \cdot m \otimes h\_2 \cdot n $$ and the coaction
\begin{equation}\label{yetter:tensor}
\delta(m\otimes n) =\sum m\_0 \otimes n\_0 \otimes n\_1m\_1.
\end{equation}

An alternative definition of the Drinfel'd double is $D'(H) =H\otimes (H^*)^{\coop}$ as coalgebras, and multiplication given by
$$(k\otimes f)(h\otimes g) = \sum kf\_1(S^{-1}(h\_1))f\_3(h\_3)h\_2\otimes f\_2g,$$
where $(\com^{\op}\otimes \id)\com^{\op}(f) = \sum f\_1\otimes f\_2 \otimes f\_3$.
We will need both definitions of the Drinfel'd double since two of the articles we will be using \cite{benkart-witherspoon, radford} use these different definitions.
The following lemma gives the relationship between these two
definitions of the Drinfel'd double.

\begin{lemma}\label{different-doubles}
$D'(H)\simeq D(H^*)^{\coop}$ as Hopf algebras.
\end{lemma}
\begin{proof}
As $H^{**}\simeq H$, we have $D(H^*)\cong H^{\coop}\otimes H^*$, with
multiplication 
\begin{eqnarray*}
(k \otimes f)(h \otimes g) &=& \sum k \left( f\_1 \rightharpoonup h \leftharpoonup (S^*)^{-1}(f\_3)\right) \otimes f\_2g \\
&=& \sum k \left( f\_1(h\_2)h\_1 \leftharpoonup (f\_3 \circ S^{-1}) \right) \otimes f\_2g \\
&=& \sum k f\_1(h\_3)(f\_3 (S^{-1}(h\_1))h\_2 \otimes f\_2g,
\end{eqnarray*}
where $(\com^{\op}\otimes \id)\com^{\op}(f)=\sum f\_3 \otimes f\_2 \otimes f\_1$.
So $D'(H) \simeq D(H^{*})$ as algebras. As coalgebras $D(H^*) \simeq H^{\coop}\otimes H^* = \left( H \otimes (H^*)^{\coop}\right)^{\coop} = D'(H)^{\coop}.$
\end{proof}

Next we describe the results from \cite{radford}. Although Radford's results are more general, we will only write them for $\k$ an algebraically closed field of characteristic 0 and only for finite-dimensional Hopf algebras.   
\begin{lemma}[Radford \cite{radford}]\label{lemma-beta-action} Let $H$ be a finite-dimensional Hopf algebra over $\k$. If $\beta \in G(H^*)$, then $H_{\beta}= (H, \db, \Delta) \in  {}_H {\mathcal YD}^H$, where 
\begin{equation}\label{beta-action}
 h\db a =\sum \langle \beta, h\_2 \rangle h\_3aS^{-1}(h\_1)
\end{equation}
for all $h,a$ in $H$.
\end{lemma}

If $\beta:H \to \k $ is an algebra map and $N$ is a right coideal of $H$, then the $H$-submodule of $H_{\beta}$ generated by N, $H\db N$, is a Yetter-Drinfel'd $H$-submodule of $H_{\beta}$. If $g\in G(H)$, then $\k g $ is a right coideal and $H\db \k g = H \db g$ is a Yetter-Drinfel'd submodule of $H_{\beta}$. For $M$ a Yetter-Drinfel'd module over $H$, $[M]$ will denote the the isomorphism class of $M$. 

\begin{prop}[Radford \cite{radford}]\label{radford} Let $H= \bigoplus_{n=0}^{\infty}H_{n}$ be a finite-dimensional graded Hopf algebra over $\k$. Suppose that $H_0=\k G$ where $G$ is a finite abelian group. Then 
$$(\beta, g) \mapsto [H\db g]$$
describes a bijective correspondence between the Cartesian product of sets $G(H^*) \times G=G(H^*)\times G(H)$  and the set of isomorphism classes of simple Yetter-Drinfel'd $H$-modules.
\end{prop}
Note that in this case, if $h \in G$, then $h \db g  = \langle \beta , h \rangle g$.

Since $H$ is finite dimensional, $H_m=H_{m+1} = \cdots = (0)$ for some $m >0$. Hence, if $\beta:H \to \k$ is an algebra map and $i>0$, then $\beta_{\mid_{H_i}} = 0$, and $\beta$ is determined by its restriction to $H_0=\k G$. 

Let \begin{equation}\label{widehatG}\widehat{G}=\Hom(G, \k^{\times}),\end{equation} the set of group homomorphisms from $G$ to $\k^{\times}= \k -\{0\}$. Then, to give an algebra map $\beta: H \to \k$, is equivalent to giving a map in $\widehat{G}$; when no confusion arises, the corresponding map in $\widehat{G}$ will also be called $\beta$.


We present a general result on the tensor product of Yetter-Drinfel'd modules. 

\begin{prop}\label{tensor products}
Let $H= \bigoplus_{n=0}^{\infty}H_{n}$ be a finite dimensional graded Hopf algebra over $\k$ with $H_0=\k G$ where $G$ is a finite abelian group.
Let $\beta , \, {\beta}' \in G(H^*)$ and $g, \, g' \in G(H)$. If $ H \db g \otimes H\DOT_{{\beta}'} g'$ is a simple Yetter-Drinfel'd module, then 
$$ H \db g \otimes H \DOT_{{\beta}'} g' \simeq H \DOT_{\beta{\beta}'} gg'$$
\end{prop}
\begin{proof}
Since $ H \db g \otimes H\DOT_{{\beta}'} g'$ is a simple Yetter-Drinfel'd module, by Proposition \ref{radford}, there exist unique ${\beta }'' \in G(H^*)$ and $g'' \in G(H)$ such that $$ H \db g \otimes H \DOT_{{\beta}'} g' \simeq H \DOT_{\beta ''} g''$$ as Yetter-Drinfel'd modules. 
Let $\Phi :  H \db g \otimes H \DOT_{{\beta}'} g' \to H\DOT_{\beta ''}g''$ be such an isomorphism. 
Since $\Phi$ is a comodule map, we have 
\begin{equation*} 
(\Phi \otimes \id )\left(\sum g\_0 \otimes g'\_0 \otimes g'\_1 g\_1\right) = \com (\Phi (g \otimes g')).
\end{equation*}
Then
\begin{equation}\label{eqn:phi}
\Phi (g \otimes g') \otimes g'g = \com (\Phi (g \otimes g')).
\end{equation}
This implies that $\k\Phi (g \otimes g')$ is a (simple) right coideal of $H\DOT_{\beta ''} g''$. In \cite{radford} it was shown that if $N$ is a simple right coideal of $H$, then the only coideal contained in $H\db N$ is $N$. 

Therefore $\k\Phi (g \otimes g') = \k g''$ and so $g'' = \lambda \Phi(g \otimes g')$ for some $0 \neq \lambda \in \k$; we may assume that $\lambda =1$. Applying $\epsilon \otimes \id $ to both sides of Equation \eqref{eqn:phi}, we get that $\Phi(g\otimes g') = \epsilon (\Phi(g \otimes g')) g'g$. We then have:
$$ g''= \Phi (g \otimes g') = \epsilon (\Phi (g \otimes g')) g'g.$$
Since distinct group-like elements are linearly independent, this implies that $g''= g'g$.

It is enough to show that $\beta''$ and $\beta \beta'$ agree on $G$.
Let $h \in G$, then
\begin{eqnarray*}
\beta '' (h) gg' &=& h \DOT_{{\beta}''} gg'  =  h \DOT_{\beta ''}(\Phi (g \otimes g')) = \Phi( h \cdot (g \otimes g'))\\
&=&  \Phi (h \db g \otimes h \DOT_{\beta '} g') = \Phi ( \beta (h) \beta '(h) g \otimes g') = (\beta  \beta ')(h) gg',
\end{eqnarray*}
and so $\beta ''(h) =  (\beta \beta ')(h)$ for all $h$ in $G$. 
\end{proof}

If $H$ is any Hopf algebra and $\gamma : H \to \k$ is an algebra map, then $\gamma$ has an inverse in $\Hom(H,\k)$ given by $\gamma^{-1}(h) = \gamma(S(h))$.

Let $N=\k n$ be a one-dimensional $H$-module. Then there is an algebra homomorphism $\gamma : H \to \k$ such that $ h \cdot n = \gamma(h) n $ for all $ h \in H$. Let $\k_{\gamma}$ be $\k$ as a vector space with the action given by $h\cdot 1 = \gamma(h)$, and so $N \simeq \k_{\gamma}$ as $H$-modules.

If $M$ is any $H$-module and $\gamma : H\to \k$ is an algebra morphism, then the natural vector space isomorphism $M \otimes \k_{\gamma} \simeq M$ endows $M$ with a new module structure, $\cdot '$, given by $h\cdot'm = \sum \gamma(h\_2)h\_1\cdot m$. We will denote this module by $M_{\gamma}.$

Note that for any $H$-module $M$, $$\left( M_{\gamma} \right)_{\gamma^{-1}}=M_{\epsilon}=M.$$
Let $\gamma :H \to \k$ be an algebra map. If $M$ is an $H$-module and $N$ is a submodule of $M$, then $N_{\gamma}$ is a submodule of $M_{\gamma}.$ In particular, $M$ is simple if and only if $M_{\gamma}$ is simple. Let $\mbox{Soc}(M)$ denote the socle of $M$, that is, $\mbox{Soc}(M)= \oplus N,$ the sum over all simple submodules of $M$. We have that $$\mbox{Soc}(M_{\gamma}) = \left(\mbox{Soc}(M)\right)_{\gamma}.$$

We give a construction that will be used to get a factorization of simple $D(H)$-modules. Given any Hopf algebra $H$ with counit $\epsilon$ and $L$ a subset of $H$, let $$L^{+}= L \cap \ker \epsilon.$$ 

Note that if $L$ is a subcoalgebra of $H$, then $L^{+}$ is a coideal and hence $H/L^{+}$ is a coalgebra. In this case, let $\langle L^+ \rangle = HL^+H$ be the two-sided ideal generated by $L^+$, then $H/\langle L^+ \rangle$ is a bialgebra. We will use this construction in the particular case when $L \subset Z(H)$, the center of $H$, in which case $\langle L^+ \rangle = HL^+ $ and so $H/HL^{+}$ is a bialgebra. If in addition $S(L^{+}) \subset L^+$, then $H/HL^+$ is a Hopf algebra.
A simple calculation shows that if $L = \k J$ with $J$ a subgroup of $G(H)$, the group of group-like elements of $H$, then $$L^+=\k\left\{g-1 : \, g \in J\right\}.$$ 
\begin{rem}\label{schneider:quotient}
In \cite{schneider} H.-J. Schneider strengthened the Nichols-Zoeller theorem and showed that if $H$ is a finite-dimensional Hopf algebra and $L$ is a Hopf subalgebra of $H$, then $H \simeq H/HL^+ \otimes L$ as right $L$-modules \cite{schneider} .
In particular $$\dim(H/HL^+) = \frac{\dim(H)}{\dim(L)}.$$
\end{rem}
\begin{definition}\label{dn:central-group-like}
For $H$ a finite-dimensional Hopf algebra, let $$G_C(H)=G(H)\cap Z(H)$$ denote the group of central group-like elements of $H$ and let $$\overline{H}= H/H(\k G_C(H))^{+}.$$
\end{definition} 
Then $\overline{H}$ is a Hopf algebra, and by Remark \ref{schneider:quotient} $$\dim(\overline{H}) = \frac{\dim(H)}{|G_C(H)|}.$$


\section{Two-parameter quantum groups}\label{quantum}

Let $\alpha_j = \epsilon_{j}- \epsilon_{j+1} \ (j = 1, \dots, n-1)$.
Let $r,s\in \k ^{\times}$ be roots of unity with $r\neq s$ and $\ell$ be the least common multiple of the orders of $r$ and $s$.
Let $\theta$ be a primitive $\ell$th root
of unity and $y$ and $z$ be nonnegative integers such that $r=\theta^y$ and $s=\theta^z$. Takeuchi defined the following Hopf algebra \cite{takeuchi}.
\begin{definition}\label{quantum-group}The algebra $U=U_{r,s}(\mathfrak{sl}_{n})$ is the unital associative
$\k$-algebra  
generated by $
\{e_j, \ f_j, \ \w_j^{\pm 1}, \ (\w_j')^{\pm 1}, \ \ 1 \leq j < n \}$, subject
to the following relations.
\medbreak
 
 \begin{itemize}
\item[(R1)]  The $\w_i^{\pm 1}, \ (\w_j')^{\pm 1}$ all commute with one
another and \\
 $\w_i \w_i^{-1}= \w_j'(\w_j')^{-1}=1,$
\item[(R2)] $ \w_i e_j = r^{\langle \epsilon_i,\alpha_j\rangle }s^{\langle
\epsilon_{i+1},\alpha_j\rangle}e_j \w_i$ \ \ and \ \ $\w_if_j =
r^{-\langle\epsilon_i,\alpha_j\rangle}s^{-\langle\epsilon_{i+1},\alpha_j\rangle} f_j\w_i,$ 
\item[(R3)]
$\w_i'e_j = r^{\langle\epsilon_{i+1},\alpha_j\rangle}s^{\langle\epsilon_{i},\alpha_j\rangle}e_j
\w_i'$ \ \ and \ \ $\w_i'f_j =
r^{-\langle\epsilon_{i+1},\alpha_j\rangle}s^{-\langle\epsilon_{i},\alpha_j\rangle}f_j \w_i'$, 
\item[(R4)]
$\displaystyle{[e_i,f_j]=\frac{\delta_{i,j}}{r-s}(\w_i-\w_i').}$
\item[(R5)] $[e_i,e_j]=[f_i,f_j]=0 \ \ \text{ if }\ \ |i-j|>1, $  
\item[(R6)]
 $e_i^2e_{i+1}-(r+s)e_ie_{i+1}e_i+rse_{i+1}e_i^2 = 0,$ \\
$e_i e^2_{i+1} -(r+s)e_{i+1}e_ie_{i+1} +rse^2_{i+1}e_i = 0,$  
\item[(R7)]
 $f_i^2f_{i+1}-(r^{-1}+s^{-1})f_if_{i+1}f_i +r^{-1}s^{-1}f_{i+1}f_i^2 =
 0,$\\
$f_i f^2_{i+1} -(r^{-1}+s^{-1})f_{i+1}f_if_{i+1}+r^{-1}s^{-1} f^2_{i+1} f_i=0, $
\end{itemize}
for all $1\leq i,j <n.$
\end{definition}  

The following coproduct, counit, and antipode give $U$ the structure of a Hopf algebra:
\begin{eqnarray*}
\Delta(e_i)=e_i\otimes 1+\omega_i\otimes e_i, &\ \ & \Delta(f_i)=1\otimes f_i +
f_i\otimes\omega_i',\\
\epsilon(e_i)=0, &\ \ & \epsilon(f_i)=0,\\
S(e_i)=-\omega_i^{-1}e_i,&\ \ & S(f_i)=-f_i(\omega_i')^{-1},
\end{eqnarray*} 
and $\omega_i,\omega_i'$ are group-like, for all $1\leq i <n$.

Let $U^0$ be the group algebra generated by all $\omega_i^{\pm 1}$,
$(\omega_i')^{\pm 1}$ and let $U^{+}$ (respectively, $U^{-}$) be the
subalgebra of $U$ generated by all $e_i$ (respectively, $f_i$). Let
$$\E_{j,j}=e_j \ \ \mbox{ and } \ \ \E_{i,j}=e_i\E_{i-1,j}-
  r^{-1}\E_{i-1,j}e_i \quad  (i>j),$$
$$\F_{j,j}=f_j \ \ \mbox{ and } \ \ \F_{i,j}=f_i\F_{i-1,j}-
s\F_{i-1,j}f_i \quad (i>j).$$
The algebra $U$ has a triangular decomposition $U\cong U^-
\otimes U^0\otimes U^+$ (as vector spaces), and the subalgebras $U^+$, $U^-$ 
respectively have monomial Poincar\'e-Birkhoff-Witt (PBW) bases
\cite{kharchenko, b-k-l}
 \begin{equation*}
\mathcal E: = \{\E_{i_1,j_1}\E_{i_2,j_2}\cdots\E_{i_p,j_p}\mid (i_1,j_1)\leq
(i_2,j_2)\leq\cdots\leq (i_p,j_p)\mbox{ lexicographically}\},
\end{equation*}
\begin{equation*}
\mathcal F:= \{\F_{i_1,j_1}\F_{i_2,j_2}\cdots\F_{i_p,j_p}\mid (i_1,j_1)\leq
(i_2,j_2)\leq\cdots\leq (i_p,j_p)\mbox{ lexicographically}\}.
\end{equation*}

It is shown in \cite{benkart-witherspoon} that all
$\E^{\ell}_{i,j}$, $\F^{\ell}_{i,j}$, $\omega_i^{\ell}-1$, and
$(\omega'_i)^{\ell}-1$ ($1\leq j \leq i< n$) are central
in $U_{r,s}({\mathfrak{sl}}_{n})$.   The ideal $I_{n}$  
generated by these elements is a Hopf ideal  
\cite[Thm.\ 2.17]{benkart-witherspoon}, and so the quotient
\begin{equation}\label{resqg}{\mathfrak{u}}={\mathfrak{u}}_{r,s}({\mathfrak{sl}}_{n})=
  U_{r,s}({\mathfrak{sl}}_{n})/I_{n}\end{equation}
is a Hopf algebra, called the {\em restricted two-parameter quantum group}.
Examination of the PBW-bases $\mathcal E$ and $\mathcal F$ shows that ${\mathfrak{u}}$ is finite-dimensional and Benkart and Witherspoon showed that ${\mathfrak{u}}$ is pointed \cite[Prop.\ 3.2]{benkart-witherspoon}.  
 \medbreak

Let $\mathcal E_\ell$ and $\F_{\ell}$ denote the sets of monomials in $\mathcal E$ and $\F$ respectively, in which each $\mathcal E_{i,j}$ or $\F_{i,j}$ appears as a factor at most $\ell-1$ times. Identifying cosets in $\mathfrak u$ with their
representatives, we may assume $\mathcal E_{\ell}$ and $\F_{\ell}$ are basis for the subalgebras of
$\mathfrak u$ generated by the elements $e_i$ and $f_i$ respectively. 

Let $\mathfrak b$ be the Hopf subalgebra of ${\mathfrak{u}}_{r,s}(
{\mathfrak{sl}}_n)$ generated by $\{\omega_i, e_i : \, 1 \leq i <n\}$, and ${\mathfrak{b}}'$ the subalgebra generated by $\{ \omega_i',\, f_i: \, 1\leq i < n\}$. 

Benkart and Witherspoon showed that, under some conditions on the parameters $r$ and $s$, ${\mathfrak{b}}^*\simeq ({\mathfrak{b}}')^{\coop}$ as Hopf algebras (\cite[Lemma\ 4.1]{benkart-witherspoon}). This implies that ${\mathfrak{b}}\simeq \left(({\mathfrak{b}}')^{\coop}\right)^*$; we present their lemma using the dual isomorphism of the original one.  
\begin{lemma}\label{action-dual}{\rm \cite[Lemma \ 4.1]{benkart-witherspoon}} If $\gcd(y^{n-1}-y^{n-2}z+\cdots +(-1)^{n-1}z^{n-1},\ell)=1$ and $rs^{-1}$ is a primitive $\ell$th root of unity, then 
${\mathfrak{b}}\simeq \left(({\mathfrak{b}}')^{\coop}\right)^*$ as Hopf algebras. Such an isomorphism is given by 
\begin{equation}\label{w-action}
\langle \w_i, \w'_j \rangle= r^{\langle \epsilon_i, \alpha_j\rangle} s^{\langle \epsilon_{i+1}, \alpha_j\rangle} \quad \mbox{and}  \quad \langle \w_i,  f_j\rangle=0,
\end{equation}
and 
\begin{equation}\label{f-action}
\langle e_i, f_j^ag \rangle = \delta_{i,j}\delta_{1,a} \quad \forall \, g \in G(\mathfrak{b}'). 
\end{equation}
\end{lemma}

\begin{prop}\label{double}
{\rm \cite[Thm.\ 4.8]{benkart-witherspoon}} Assume $r = \theta^y$ and $s = \theta^z$,  
where $\theta$ is a primitive $\ell$th root of unity, and 
\begin{equation*}\label{relprime} \gcd(y^{n-1}-y^{n-2}z+\cdots +(-1)^{n-1}z^{n-1},\ell)=1.
\end{equation*}  
Then there is an
isomorphism of Hopf algebras ${\mathfrak{u}}_{r,s}({\mathfrak{sl}}_{n}) 
\cong D'({\mathfrak{b}}).$
\end{prop} 
   
Under the assumption that $\gcd(y^{n-1}-y^{n-2}z+\cdots +(-1)^{n-1}z^{n-1},\ell)=1$, combining lemmas \ref{different-doubles} and \ref{w-action}, and Proposition \ref{double}, we get that $\ursn \simeq (D((\mathfrak{b}')^{\coop}))^{\coop}$. Hence, $\ursn$-modules are Yetter-Drinfel'd modules for $(\mathfrak{b}')^{\coop}$ (only the algebra structure of $\ursn$ plays a role when studying $\ursn$-modules, hence $\ursn$-modules are $D((\mathfrak{b}')^{\coop})$-modules). For simplicity we will denote $H=(\mathfrak{b}')^{\coop}$. Then $G=G(H)=\langle \w_i' : \, 1 \leq i < n\rangle$ and $H$ is a graded Hopf algebra with $\w_i' \in H_0$, and $f_i \in H_1$ for all $1 \leq i <n$. Therefore Proposition \ref{radford} applies to $H$ and isomorphism classes of $\ursn$-modules (or simple Yetter-Drinfel'd $H$-modules) are in one to one correspondence with $G(H^*) \times G(H)$.

\subsection[Factorization of simple $\ursn$-modules]{Factorization of simple $\ursn$-modules}
We study under what conditions a simple $\ursn$-module can be factored as the tensor product of a one-dimensional module and a simple module which is also a module for $\overline{\ursn}= \ursn/\ursn(\k G_C(\ursn))^+.$ Let $\ell$, $n$, $y$ and $z$ be fixed and $\theta$ be a primitive $\ell$th root of unity. Let $A$ be the $(n-1) \times (n-1) $ matrix 
\begin{equation*}
A 
= 
\left( 
\begin{array}{cccccc}
y-z &      z     & 0 & 0& \cdots & 0 \\
-y & y-z & z & 0 &  \cdots & 0 \\
\vdots &&&&& \vdots \\ 
0&  \cdots& 0 & -y & y-z & z \\
0& \cdots & \cdots &0& -y & y-z
\end{array}
\right)
\end{equation*}

The determinant of $A$ is $y^{n-1} -  y^{n-2}z + \cdots + (-1)^{n-1}z^{n-1}$. Throughout this section, assume that $\gcd(y^{n-1} -  y^{n-2}z + \cdots + (-1)^{n-1}z^{n-1}, \ell)=1$, and so $\mbox{det}(A)$ is invertible in $\Z / \ell \Z$.
We start by describing the set of central group-like elements in $\ursn$. Clearly $G(\ursn)=\langle \w_i, \w_i' : 1 \leq i <n\rangle$.

\begin{prop}\label{centralu} A group-like element $g=\w_1^{a_1} \cdots\w_{n-1}^{a_{n-1}} \w_1'^{b_1} \cdots \w_{n-1}'^{b_{n-1}}$ is central in $\ursn$ if and only if 
\begin{equation*}
\left( \begin{array}{c} b_1 \\ \vdots \\ b_{n-1} \end{array} \right)
=
A^{-1}A^{\mbox{t}}\left( \begin{array}{c}a_1 \\ \vdots \\ a_{n-1} \end{array}\right)
\end{equation*}
in $\left( \Z/\ell\Z \right)^{n-1}$.
\end{prop}
 
\begin{proof}
The element $g$ is central in $\ursn$ if and only if $ge_k = e_kg$ and $gf_k = f_kg$ for all $k=1, \cdots, n-1$. 
By the relations (R2) and (R3) of the definition of $U_{r,s}(\mathfrak{sl}_n)$, for all $k=1, \cdots, n-1$ we have that
\begin{eqnarray*}
ge_k &=& \prod_{i=1}^{n-1} \left(r^{\langle \epsilon_i, \alpha_k \rangle} s^{\langle \epsilon_{i+1}, \alpha_k \rangle}\right) ^{a_i} \prod_{j=1}^{n-1} \left( r^{\langle \epsilon_{j+1}, \alpha_k \rangle} s^{\langle \epsilon_{j}, \alpha_k \rangle}\right)^{b_j} e_kg \mbox{ and} \\
gf_k &=& \prod_{i=1}^{n-1} \left(r^{-\langle \epsilon_i, \alpha_k \rangle} s^{-\langle \epsilon_{i+1}, \alpha_k \rangle}\right)^{a_i} \prod_{j=1}^{n-1} \left(r^{-\langle \epsilon_{j+1}, \alpha_k \rangle} s^{-\langle \epsilon_{j}, \alpha_k \rangle}\right)^{b_j}f_kg.
\end{eqnarray*}
Then $g$ is central if and only if
\begin{eqnarray*}
1 &=&  \prod_{i=1}^{n-1} \left(r^{\langle \epsilon_i, \alpha_k \rangle} s^{\langle \epsilon_{i+1}, \alpha_k \rangle}\right)^{a_i} \prod_{j=1}^{n-1} \left( r^{\langle \epsilon_{j+1}, \alpha_k \rangle} s^{\langle \epsilon_{j}, \alpha_k \rangle} \right)^{b_j} \\
&=&
s^{a_{k-1}}r^{a_k}s^{-a_k}r^{-a_{k+1}}r^{b_{k-1}}r^{-b_k}s^{b_k}s^{-b_{k+1}}, \, \forall k=1,\cdots, n-1, \end{eqnarray*} where $a_0=a_n=0=b_0=b_n$.
Since $r= \theta ^y$ and $s= \theta ^z$, the last equation holds if and only if
\begin{equation}\label{cond}z a_{k-1} + (y-z)a_k - y a_{k+1}=\left(-yb_{k-1} + (y-z)b_k + zb_{k+1}\right)\mbox{ mod}\ell, \end{equation}
for all $k=1, \cdots, n-1$; 
that is, if and only if 
\begin{equation*} A^{\mbox{t}} \left( \begin{array}{c} a_1 \\ \vdots \\ a_{n-1} \end{array} \right)
=
A \left( \begin{array}{c} b_1 \\ \vdots \\ b_{n-1} \end{array} \right)
\end{equation*} 
in $\left( \Z / \ell \Z \right)^{n-1}.$
\end{proof}

\begin{ex}
For $\mathfrak{u}_{\theta, \theta^{-1}}(\mathfrak{sl}_n)$ ($y=1$ and $z=\ell-1$), the matrix $A$ is symmetric.
Therefore, a group-like element $g=\w_1^{a_1} \cdots\w_{n-1}^{a_{n-1}} \w_1'^{b_1} \cdots \w_{n-1}'^{b_{n-1}}$ is central if and only if $b_i=a_i$ for all $i=1, \cdots, n-1$.
\end{ex}
Recall that $\ursn (\k G_C(\ursn))^+ = \ursn \{ g - 1 : \, g \in G_C(\ursn) \} $.
In particular, by the last example, we have that $\mathfrak{u}_{\theta,\theta^{-1}} (\k G_C(\mathfrak{u}_{\theta,\theta^{-1}}(\mathfrak{sl}_{n})))^+ $ is generated by $\left\{ \w_i^{-1} - \w'_i : \, i=1, \ldots, n-1 \right\}$. This gives $\overline{\mathfrak{u}_{\theta,\theta^{-1}}} \simeq \mathfrak{u}_{\theta}(\mathfrak{sl}_n)$, the one parameter quantum group.

Henceforth $r$ and $s$ are such that $rs^{-1}$ is also a primitive $\ell$th root of unity, that is, $\gcd(y-z, \ell)=1$. 
\begin{rem} If $\beta \in G(H^*)$ and $g = \w_1'^{c_1}  \cdots \w_{n-1}'^{c_{n-1}}  \in G(H)$, by propositions \ref{yetter} and \ref{double}, the Yetter-Drinfel'd module $H\db g$ is also a $\ursn$-module where the action of $\mathfrak{b}$ is given by combining equations \eqref{action-dual-general} and \eqref{w-action}.In particular,
$$\w_i \cdot g = \langle \w_i, g \rangle g = \prod_{j=1}^{n-1}{ \langle \w_i , \w_j' \rangle}^{c_j}g = \prod_{j=1}^{n-1}\left( r^{\langle\epsilon_i, \alpha_j\rangle}s^{\langle\epsilon_{i+1}, \alpha_j \rangle} \right)^{c_j}g.
$$
\end{rem}


\begin{prop}
Let $\beta \in G(H^*)$ be defined by $\beta({\w_i'}) = \theta ^{ \beta_i}$ and $g= \w_1'^{c_1}  \cdots \w_{n-1}'^{c_{n-1}}$. The simple $\ursn$-module $H\db g$ is naturally a $\overline{\ursn}$-module if and only if  
\begin{equation}\label{betas}
\left( \begin{array}{c}
\beta_1 \\ \vdots \\ \beta_{n-1} \end{array} \right)
=
-A^{\mbox{t}} \left( \begin{array}{c} c_1 \\ \vdots \\c_{n-1} \end{array} \right) 
\end{equation}
in $(\Z/ \ell Z)^{n-1}.$
\end{prop} 

\begin{proof}
$H\db g$ is a $\overline{\ursn}$-module if and only if $(h-1) \cdot m = 0$ for all $h$ in $G_C(\ursn)$ and $m$ in $H\db g$. If $h \in G_C(\ursn)$, then $h \cdot m =m$ for all $m$ in $H \db g$ if and only if $h \cdot g = g$. 

Let $h=\w_1^{a_1}\cdots \w_{n-1}^{a_{n-1}}\w_1'^{b_1}\cdots \w_{n-1}'^{b_{n-1}} \in G_C(\ursn)$; then by Proposition \ref{centralu}  $$
\left( \begin{array}{c} b_1 \\ \vdots \\ b_{n-1} \end{array} \right) 
=
A^{-1}A^t\left( \begin{array}{c} a_1 \\  \vdots \\a_{n-1} \end{array} \right).$$

We have
\begin{equation}\label{eqw'}
\begin{split}
\w_1'^{b_1}\cdots \w_{n-1}'^{b_{n-1}}\db g &= \beta(\w_1'^{b_1}\cdots \w_{n-1}'^{b_{n-1}})g \\
 &= \theta^{b_1\beta_1 + \cdots + b_{n-1}\beta_{n-1}}g \end{split}\end{equation}
and
\begin{equation}\label{eqw}
\begin{split}
\w_1^{a_1}\cdots \w_{n-1}^{a_{n-1}} \cdot g &= \langle \w_1^{a_1}\cdots \w_{n-1}^{a_{n-1}}, g \rangle g \\
&= \prod_{i=1}^{n-1} {\langle \w_i , g \rangle}^{a_i}g \\
&= \prod_{i=1}^{n-1}\prod_{j=1}^{n-1}\left( r^{\langle\epsilon_i, \alpha_j\rangle}s^{\langle\epsilon_{i+1}, \alpha_j \rangle} \right)^{a_ic_j}g \\
&=  \prod_{i=1}^{n-1}\left( r^{-c_{i-1}}\left(rs^{-1}\right)^{c_i}s^{c_{i+1}} \right)^{a_i}  g \\
&= \theta^xg
\end{split}
\end{equation} where $c_0=c_n=0$ and $x=\sum_{i=1}^{n-1} \left( -yc_{i-1} + (y-z)c_i + zc_{i+1} \right)a_i$.
From equations \eqref{eqw'} and \eqref{eqw} we get that $$h \cdot g = \theta^{x + \sum_{i=1}^{n-1}b_i\beta_i}g$$ 
for all $\left( \begin{array}{c} a_1 \\  \vdots \\ a_{n-1} \end{array} \right), \mbox{ where }
\left( \begin{array}{c} b_1 \\ \vdots \\ b_{n-1} \end{array} \right) 
=
A^{-1}A^t\left( \begin{array}{c} a_1 \\  \vdots \\a_{n-1} \end{array} \right)$.
Now $$\sum_{i=1}^{n-1} \left( -yc_{i-1} + (y-z)c_i + zc_{i+1} \right)a_i+ \sum_{i=1}^{n-1}b_i\beta_i= 0 \mbox{ mod} \ell$$ if and only if
\begin{eqnarray*}
{\left( \begin{array}{c} a_1 \\  \vdots \\a_{n-1} \end{array} \right)}^{\mbox{t}} A \left( \begin{array}{c} c_1 \\ \vdots \\ c_{n-1} \end{array} \right) &=& - \left( \begin{array}{c} b_1 \\ \vdots \\ b_{n-1} \end{array} \right)^{\mbox{t}} \left( \begin{array}{c} \beta_1 \\ \vdots \\ \beta_{n-1} \end{array} \right) \mbox{ in } \Z/\ell \Z . \\
\end{eqnarray*}
We then have that $H\db g$ is a $\overline{\ursn}$-module, if and only if
\begin{eqnarray*}
{\left( \begin{array}{c} a_1 \\  \vdots \\a_{n-1} \end{array} \right)}^{\mbox{t}} A \left( \begin{array}{c} c_1 \\ \vdots \\ c_{n-1} \end{array} \right) &=& -\left( \begin{array}{c} a_1 \\ \vdots \\ a_{n-1} \end{array} \right)^{\mbox{t}} A \left( A^{-1} \right)^{\mbox{t}} \left( \begin{array}{c} \beta_1 \\ \vdots \\ \beta_{n-1} \end{array} \right),
\end{eqnarray*}
for all $(a_1, \cdots, a_{n-1})$ in $\left( \Z / \ell\Z \right)^{n-1}$. This occurs if and only if
\begin{eqnarray*} 
\left( \begin{array}{c} c_1 \\ \vdots \\ c_{n-1}  \end{array} \right) &=&-\left( A^{-1} \right)^{\mbox{t}} \left( \begin{array}{c} \beta_1 \\ \vdots \\ \beta_{n-1} \end{array} \right) \mbox{in } (\Z/ \ell \Z)^{n-1}.
\end{eqnarray*}
\end{proof}
Given $g =(\w'_1)^{c_1} \cdots (\w'_{n-1})^{c_{n-1}}\in G(H)$, let $\beta_1, \dots, \beta_n$ be defined as in Equation \eqref{betas}. We will denote by $\beta_g$ the algebra map given by $\beta_g(\w_i')= \theta^{\beta_i}$.

For any Hopf algebra $H$, let $\mathcal{S}_H$ denote the denote the set of isomorphism classes of simple $H$-modules. Then  $\mathcal{S}_{\overline{H}}$ can be identified as the subset of $\mathcal{S}_H$ consisting of the $H$-modules that are naturally $\overline{H}$-modules. Combining the last proposition with Proposition \ref{radford}, we get
\begin{cor}\label{quotient}
The correspondence $G(H) \to \mathcal{S}_{\baru}$ given by
$$g \mapsto [H \DOT_{\beta_g}g]$$
is a bijection.
\end{cor}

\begin{ex}\label{u2} In the $\mathfrak{u}_{\theta,\theta^{-1}}(\mathfrak{sl}_2)$ case, the matrix $A$ is $A=(2)$. Then, the simple  $\mathfrak{u}_{\theta,\theta^{-1}}(\mathfrak{sl}_2)$-modules that are naturally $\mathfrak{u}_{\theta}(\mathfrak{sl}_2)$-modules, are of the form $H\db (\w')^c$ with $\beta(\w')={\theta}^{-2c}$.
\end{ex}

Next we describe the 1-dimensional $\ursn$-modules in terms of Radford's construction. For an algebra map $\chi:\ursn \to \k$, let $\k_{\chi}$ be the 1-dimensional $\ursn$-module given by $h \cdot 1 = \chi(h)1$. Since $e_i^{\ell}=0=f_i^{\ell}$ we have that $\chi(e_i) = \chi(f_i)=0$, and this together with (R4) of the Definition \ref{quantum-group} of $U_{r,s}(\mathfrak{sl}_n)$, gives that $\chi(\w_i) = \chi(\w_i')$. For each $i=1, \ldots, n-1$, since $\w_i^{\ell}=1$, $\chi(\w_i)=\theta^{\chi_i}$ for some $0 \leq \chi_i < \ell$.
\begin{prop}\label{1dim-prop}
If $\chi:\ursn \to \k$ is an algebra map, then $\k_{\chi} \simeq H \DOT_{\chi_{\mid_{H}}} g_{\chi}$, where $$g_{\chi}= \w_1'^{d_1} \cdots \w_{n-1}'^{d_{n-1}}, \mbox{ with }$$ 
\begin{equation*}
\left( \begin{array}{c} d_1 \\ \vdots \\ d_{n-1} \end{array} \right) = A^{-1} \left( \begin{array}{c}\chi_1 \\ \vdots \\ \chi_{n-1} \end{array} \right) \mbox{ in } (\Z /\ell \Z )^{n-1} .
\end{equation*}
\end{prop}

\begin{proof}
Since $\k_{\chi}$ is a simple $\ursn$-module, we have that $\k_{\chi} \simeq H \db g$  for some unique $\beta \in G(H^*)$ and $g \in G(H)$. Let $\phi : \k_{\chi} \to H \db g$ be an isomorphism of Yetter-Drinfel'd modules. We may assume that $g=\phi(1)$; then
\begin{equation*}
\beta(\w_i')g = \w_i' \db g = \w_i' \db (\phi(1)) = \phi \left( \w_i' \cdot 1 \right) =  \phi (\chi(\w_i')1) = \chi(\w_i') g.
\end{equation*}
Therefore $\beta(\w_i')= \chi(\w_i')$ and since $\beta(f_i)=0 = \chi(f_i)$ for all $i=1, \cdots, n-1$, we have $\beta = \chi_{\mid_H}$.

We have that
\begin{equation}\label{waction1}
\begin{split}
\w_i \cdot g &= \langle \w_i, g \rangle g \\
&= \left( \prod_{j=1}^{n-1} \langle \w_i , \w_i' \rangle ^{d_j} \right) g \\
&= \left( \prod_{j=1}^{n-1} \left( r^{\langle \epsilon_i, \alpha_j \rangle}s^{\langle \epsilon_{i+1}, \alpha_j \rangle} \right)^{d_j} \right) g \\ 
&= r^{-d_{i-1}}(rs^{-1})^{d_i}s^{d_{i+1}}g \\
&= \theta^{y(d_i-d_{i-1}) + z(d_{i+1}-d_i)}g.
\end{split}
\end{equation}

On the other hand
\begin{equation}\label{waction2}
\w_i \cdot g = \w_i \cdot \phi(1) =  \phi( \w_i \cdot 1)  =  \phi (\theta^{\chi_i} 1) = \theta^{\chi_i} g.
\end{equation}
By equations \eqref{waction1} and \eqref{waction2} we have that
\begin{eqnarray*}
-y d_{i-1} + (y-z) d_i + zd_{i+1} &=& \chi_i \mbox{ mod} \ell , \; \forall i=1, \cdots , n-1; \mbox{ and so}\\
A \left( \begin{array}{c}d_1 \\ \vdots \\ d_{n-1} \end{array} \right)
&=&
\left( \begin{array}{c} \chi_1 \\ \vdots \\ \chi_{n-1} \end{array} \right) \mbox{ in } (\Z / \ell\Z)^{n-1}.
\end{eqnarray*}
\end{proof}

For any Hopf algebra $H$, let $\mathcal{S}_{H}^1= \{[ N] \in \mathcal{S}_{H} \ : \, \dim(N)=1 \}$. Combining the last proposition and Proposition \ref{radford} we get
\begin{cor}\label{1dim}
The correspondence $G(\ursn^*) \to \mathcal{S}_{\ursn}^1$ given by
$$\chi \mapsto [H\DOT_{\chi_{\mid_{H}}}g_{\chi}]$$ is a bijection.
\end{cor}
\begin{thm}\label{final-ursn}
The map $\Phi:\mathcal{S}_{\baru}  \times \mathcal{S}_{\ursn}^1 \to \mathcal{S}_{\ursn}$ given by
$$\Phi([M],[N]) =[M\otimes N] $$ is a bijection if and only if $\gcd((y-z)n, \ell)=1$.
\end{thm} 
\begin{proof}
By the last corollary we have that 1-dimensional simple $\ursn$-modules are of the form $H\DOT_{\chi_{\mid_H}}g_{\chi}$ with $\chi \in G(\ursn^*)$. Also by Corollary \ref{quotient}, simple $\baru$-modules are of the form $H \DOT_{\beta_g}g$ for $g \in G(H)$. Furthermore by Proposition \ref{tensor products}, we have that $H \DOT_{\beta_g}g \otimes H\DOT_{\chi_{\mid_H}}g_{\chi} \simeq H\DOT_{\beta_g \chi_{\mid_H}}gg_{\chi}$. Then $\Phi$ is a bijection if and only if
$$ \Psi: \left\{ (g,\beta_g): g \in G(H) \right\} \times \left\{ (g_{\chi}, \chi): \chi \in G(\ursn^*) \right\} \to G(H) \times G(H^*) $$ given by $\Psi \left( (g, \beta_g), (g_{\chi}, \chi) \right) = (gg_{\chi}, \beta_g  \chi_{\mid_{H}} )$ is a bijection.
The latter holds if and only if for all $ h=\w_1'^{b_1} \cdots \w_{n-1}'^{b_{n-1}}$ and $\gamma$ given by $\gamma({\w_i'})= \theta^{\gamma_i}$, there exist unique $g= \w_1'^{c_1} \cdots \w_{n-1}'^{c_{n-1}} $ and $\chi$ with $\chi(w_i)=\chi(\w_i') = \theta^{\chi_i}$, so that 
$h=gg_{\chi}$ and $\gamma = \beta_g \chi_{\mid_{H}}$. If $\beta_g(\w_i') = \theta^{\beta_i}$ and $g_{\chi}=\w_1'^{d_1} \cdots \w_{n-1}'^{d_{n-1}}$, then
$$
gg_{\chi} = \w_1'^{c_1} \cdots \w_{n-1}'^{c_{n-1}}\w_1'^{d_1} \cdots \w_{n-1}'^{d_{n-1}}\mbox{ and }( \beta_g  \chi_{\mid_H})(\w_i')= \theta^{\beta_i + \chi_i}.$$
Then $\Psi$ is bijective if and only if the system of equations

\begin{eqnarray*}
\left( \begin{array}{c} c_1 + d_1 \\ \vdots\\ c_{n-1} + d_{n-1} \end{array} \right) &=& \left( \begin{array}{c} b_1 \\ \vdots \\ b_{n-1}  \end{array} \right) \\
\left( \begin{array}{c}\beta_1 + \chi_1 \\ \vdots \\ \beta_{n-1} + \chi_{n-1} \end{array} \right) &=& \left( \begin{array}{c} \gamma_1 \\ \vdots \\ \gamma_{n-1} \end{array} \right)
\end{eqnarray*}
subject to
\begin{eqnarray*}
\left( \begin{array}{c} d_1 \\ \vdots \\ d_{n-1} \end{array} \right) &=& A^{-1} \left( \begin{array}{c}\chi_1 \\ \vdots \\ \chi_{n-1} \end{array} \right) \\
\left( \begin{array}{c}
\beta_1 \\ \vdots \\ \beta_{n-1} \end{array} \right)
&=&-A^{\mbox{t}} \left( \begin{array}{c} c_1 \\ \vdots \\c_{n-1} \end{array} \right)
\end{eqnarray*}
has a unique solution for all $(b_1. \cdots , b_{n-1})$, $(\gamma_1, \cdots \gamma_{n-1})$.
The last four vector equations are equivalent to 
\begin{eqnarray*}
  \left(\begin{array}{c} c_1 \\ \vdots  \\ c_{n-1} \end{array} \right)+ A^{-1} \left( \begin{array}{c} \chi_1 \\ \vdots \\ \chi_{n-1} \end{array} \right) &=&\left( \begin{array}{c} b_1 \\ \vdots \\ b_{n-1} \end{array} \right) \\
 -A^{\mbox{t}} \left( \begin{array}{c} c_1 \\ \vdots \\ c_{n-1} \end{array} \right) + \left( \begin{array}{c} \chi_1 \\ \vdots \\\chi_{n-1} \end{array} \right)&=& \left( \begin{array}{c} \gamma_1 \\ \vdots \\\gamma_{n-1} \end{array} \right)  \end{eqnarray*}
which can be written as
\begin{equation*}
\left( \begin{array}{lr} \Id & A^{-1} \\ -A^{\mbox{t}} & \Id \end{array} \right) \left( \begin{array}{c}
c_1 \\ \vdots \\ c_{n-1} \\ \chi_1 \\ \vdots\\ \chi_{n-1} \end{array} \right) = \left( \begin{array}{c} b_1 \\ \vdots\\ b_{n-1} \\ \gamma_1 \\ \vdots \\ \gamma_{n-1} \end{array} \right).
\end{equation*}
This last system has a unique solution if and only if the matrix 
\[
M=
\left( \begin{array}{lr} \Id & A^{-1} \\ -A^{\mbox{t}} & \Id \end{array} \right) 
\] is invertible in $\mbox{M}_{(n-1) \times(n-1)}\left(\Z / \ell\Z \right)$, or equivalently, if $\gcd(\mbox{det}(M), \ell)=1$. By row-reducing $M$ we have that

\begin{eqnarray*}
\mbox{det}\left( \begin{array}{lr} \Id & A^{-1} \\ -A^{\mbox{t}} & \Id \end{array} \right) &=&
\mbox{det}\left( \begin{array}{lr} A & \Id \\ -A^{\mbox{t}}& \Id \end{array} \right)\\ &=&
\mbox{det}\left( \begin{array}{lr} A+A^{\mbox{t}} & 0 \\ -A^{\mbox{t}}& \Id \end{array} \right)\\&=& \mbox{det}(A + A^{\mbox{t}}).
\end{eqnarray*}
Now
\begin{eqnarray*}
A+A^{\mbox{t}}&=&\left( 
\begin{array}{cccccc}
2(y-z) &      z-y     & 0 & 0& \cdots & 0 \\
z-y & 2(y-z) & z-y & 0 &  \cdots & 0 \\
0 & z-y & 2(y-z) & z-y & \cdots & 0\\
\vdots &&&&& \vdots \\
0&  \cdots& 0 & z-y & 2(y-z) & z-y \\
0& \cdots & \cdots &0& z-y & 2(y-z)
\end{array}
\right)
\end{eqnarray*}
\begin{eqnarray*}
&=& (y-z)\left( 
\begin{array}{cccccc}
2 &     -1     & 0 & 0& \cdots & 0 \\
-1 & 2 & -1 & 0 &  \cdots & 0 \\
0  & -1&  2 & -1&  \dots &0\\
\vdots &&&&& \vdots \\
0&  \cdots& 0 & -1 & 2 & -1 \\
0& \cdots & \cdots &0& -1 & 2
\end{array}
\right).
\end{eqnarray*}
Therefore $\mbox{det}(A+A^{\mbox{t}})= (y-z)^{n-1}n$.
We then have that $\Phi$ is a bijection if and only if $\gcd\left( (y-z)n, \ell \right)=1$.
\end{proof}
This factorization was used in \cite{pereira_thesis} to construct $\mathfrak{u}_{r,s}(\mathfrak{sl}_3)$-modules for different values of the parameters, and to conjecture a general formula for their dimensions.

\section{Pointed Hopf algebras of rank one}\label{pointed}
Recently Andruskiewitsch and Schneider classified the finite-dimensional pointed Hopf algebras with abelian groups of group-like elements, over an algebraically closed field of characteristic 0 \cite{a-s-1}. Earlier, in 2005, Krop and Radford classified the pointed Hopf algebras of rank one, where $\mbox{rank}(H)+1$ is the rank of $H\_1$ as an $H\_0$-module and $H$ is generated by $H\_1$ as an algebra, where $H\_1$ is the first term of the coradical filtration of $H$ \cite{krop-radford}. They also studied the representation theory of $D(H)$ in a fundamental case. Using Radford's construction of simple modules, in Theorem \ref{last-pointed}, we give necessary and sufficient conditions for the tensor product of two $D(H)$-modules to be completely reducible. 

\subsection{Pointed Hopf algebras of rank one of nilpotent type}

Let $G$ be a finite abelian group, $\chi: G \to \k$ a character and $ a \in G$; we call the triple $\D=(G,\chi, a)$ {\em data}.
Let $\ell:=|\chi(a)|$, $N:=|a|$ and $M= |\chi|$; note that $\ell$ divides both $N$ and $M$. In \cite{krop-radford} Krop and Radford defined the following Hopf algebra.
\begin{definition}\label{definition:HD}
Let $\D=(G,\chi, a)$ be data. The Hopf algebra $H_{\D}$ is generated by $G$ and $x$ as a $\k$-algebra, with relations:
\begin{enumerate}
\item $x^{\ell}=0$.
\item $xg= \chi(g) gx$, for all $ g \in G$.
\end{enumerate}
The coalgebra structure is given by $\com(x) =x \otimes a + 1 \otimes x$ and $\com(g) = g \otimes g$ for all $g \in G$.
\end{definition}
The Hopf algebra $H_{\D}$ is pointed of rank one. 
Let $\Gamma = \Hom(G, \k^{\times})$, the set of group homomorphisms from $G$ to $\k^{\times}$ also written $\widehat{G}$.  
\begin{prop}[Krop and Radford \cite{krop-radford}]\label{definition:dualHD}
As a $\k$-algebra, $H_{\D}^*$ is generated by $\Gamma$ and $\xi$ subject to relations:
\begin{enumerate}
\item  $\xi ^{\ell}=0$.
\item  $ \xi \gamma = \gamma(a) \gamma \xi $, for all $ \gamma \in \Gamma$.
\end{enumerate}
The coalgebra structure of $H_{\D}^*$ is determined by $\com(\xi) = \xi \otimes \chi + 1 \otimes \xi$ and $ \com(\gamma) = \gamma \otimes \gamma$ for all $ \gamma \in \Gamma$.
\end{prop}

\begin{prop}[Krop and Radford \cite{krop-radford}]\label{double-HD}
The double $D(H_{\D})$ is generated by $G, \; x, \; \Gamma, \; \xi$ subject to the relations defining $H_{\D}$ and $H_{\D}^*$ and the following relations:
\begin{enumerate}\item $g\gamma= \gamma  g $ for all $g \in G$ and $\gamma \in \Gamma$.
\item $\xi g = \chi^{-1}(g) g \xi$ for all $ g \in G$.
\item $[x, \xi] = a - \chi$.
\item $ \gamma(a) x \gamma = \gamma x$ for all $ \gamma \in \Gamma$.
\end{enumerate}
\end{prop}
Recall that the coalgebra structure of $H_{\D}^*$ in $D(H_{\D})$ is the co-opposite to the one in $H^*$. Then in $D(H_{\D})$, $\com(\xi)= \chi \otimes \xi + \xi \otimes 1$.
Note that $H_{\D}$ satisfies the hypothesis of Proposition \ref{radford}, where elements in $G$ have degree 0 and $x$ has degree 1. Therefore, simple $D(H_{\D})$-modules are of the form $H\db g$, for $g \in G$ and $\beta \in G(H^*)=\Gamma$. 

\subsection{Factorization of simple $D(H_{\mathcal{D}})$-modules}
We study under what conditions a simple $D(H_{\D})$-module can be factored as the tensor product of a one-dimensional module with a simple module which is also a module for $\overline{D(H_{\D})}=D(H_{\D})/D(H_{\D})(\k G_C(D(H_{\D})))^+$. We also study, under certain conditions on the parameters, the reducibility of the tensor product of two simple $D(H_{\D})$-modules.
 
We start by describing the central group-like elements of $D(H_{\D})$. It is clear that $G(D(H_{\mathcal{D}}))= G \times \Gamma$. An element $(g,\gamma) \in G \times \Gamma$ will be denoted by $g\gamma$.
An element $g\gamma$ is central in $D(H_{\mathcal{D}})$ if and only if
$(g\gamma) x = x(g\gamma)$ and $ (g\gamma)\xi = \xi (g\gamma)$. Using the relations of $D(H_{D})$, we have that$$g \gamma x = \gamma(a) g x \gamma = \chi^{-1}(g)\gamma(a) x g \gamma,$$ and
$$ g \gamma \xi = \gamma g \xi  = \chi(g) \gamma \xi g = \chi(g) \gamma(a)^{-1} \xi \gamma g.$$ Hence, $g\gamma$ is central if  only if $\chi^{-1}(g)\gamma(a) =1$. Let $\mbox{ev}_{\chi^{-1}a}: G \times \Gamma \to \k^{\times}$ be the character given by $\mbox{ev}_{\chi^{-1}a}(g\gamma)=\chi^{-1}(g)\gamma(a)$; we just showed the following lemma:
\begin{lemma}\label{centrals}
 $G_C(D(H_{\mathcal{D}})) = \ker(\mbox{ev}_{\chi^{-1}a})$.
\end{lemma}

Let $\alpha:D(H_{\mathcal{D}}) \to \k$ be an algebra map; then $\alpha(x)=\alpha(\eta)=0$ (because $0=x^{\ell}=\xi^{\ell}$) and $\alpha(a)=\alpha(\chi)$ (by the third relation in Proposition \ref{double-HD}). Since  $\alpha(x)=\alpha(\eta)=0$, we can think of $\alpha$ as a group homomorphism $\alpha: G \times \Gamma \to \k^{\times}$, that is, $\alpha \in \widehat{G\times \Gamma}\simeq\Gamma\times G$. Let $\beta_{\alpha} \in \Gamma$ and $g_{\alpha} \in G$ be such that $\alpha=\beta_{\alpha}g_{\alpha}$; that is $\alpha(g\gamma)=\beta_{\alpha}(g)\gamma(g_{\alpha})$ for all $g\gamma $ in $G \times \Gamma$. If we extend $\beta_{\alpha}$ to $H_{\D}$ by setting $\beta_{\alpha}(x) =0$ and also call this extension $\beta_{\alpha}$ (as no confusion will arise), we have $\beta_{\alpha}= \alpha_{|_{H_{\D}}}$.
Let $\k_{\alpha}$ be the one-dimensional module defined by $h\cdot k = \alpha(h)k$ for all $h \in D(H_{\mathcal{D}})$ and $k \in \k$.

\begin{prop}\label{1dimK-R}
$\k_{\alpha} \simeq H_{\mathcal D} \DOT_{\beta_{\alpha}} g_{\alpha}$ as Yetter-Drinfel'd $H_{\D}$-modules.
\end{prop}
\begin{proof}
Since $\k_{\alpha}$ is a simple Yetter-Drinfel'd module, there exists an isomorphism of Yetter-Drinfel'd modules $\Phi : \k_{\alpha} \to H_{\mathcal D} \db g $ for some algebra map $\beta: H_{\mathcal D} \to \k$ and some $g \in G$.
We may assume that $\Phi(1)=g$. Let $h \in G$, we have 
$$h \db g = \beta(h) g.$$ 
Since $ \Phi$ is a module map,
\begin{eqnarray*}
h\db g &=& h \db \Phi(1) = \Phi (h \cdot 1) = \Phi( \alpha(h))\\
&=& \alpha(h) \Phi(1) = \beta_{\alpha}(h)g.
\end{eqnarray*}
We then have $\beta(h) = \beta_{\alpha}(h)$ for all $h$ in $G$, and since $\beta(x) = \beta_{\alpha}(x)=0$, $\beta= \beta_{\alpha}$.

If $\gamma \in \Gamma$, then
\begin{equation*}
\gamma \db g = \gamma( g)  g.
\end{equation*}
On the other hand,
\begin{equation*}
\gamma \db g = \gamma \db \Phi(1) = \Phi (\gamma \cdot 1) = \Phi( \alpha(\gamma) 1)
= \alpha(\gamma)\Phi(1) = \gamma (g_\alpha)g.
\end{equation*}
Then $\gamma(g)= \gamma(g_{\alpha})$ for all $\gamma \in \Gamma$, hence $g=g_{\alpha}$.
\end{proof}
For simplicity let $K=\ker(ev_{\chi^{-1}a})$.
If $\alpha= \beta_{\alpha} g_{\alpha} \in \Gamma \times G$, the condition $\alpha(a) = \alpha(\chi)$ is $\beta_{\alpha}(a)=\chi(g_{\alpha})$ or $\chi^{-1}(g_{\alpha})\beta_{\alpha}(a) =1$. Hence, $\alpha$ in $\Gamma \times G$  defines a one-dimensional module if and only if $g_{\alpha} \beta_{\alpha} \in \ker(\mbox{\rm ev}_{\chi^{-1}a})=K$. This, together with the previous proposition, shows 

\begin{cor}\label{1-dim-KR}
The set $\mathcal{S}_{D(H_{\D})}^1$ of isomorphism classes of one dimensional $D(H_{\mathcal D})$-modules is in one to one correspondence with $K$. 
\end{cor}

Recall that $\overline{D(H_{\D})}=D(H_{\D})/D(H_{\D})(\k G_C(D(H_{\D})))^+$. Since $G_C(D(H_{\D}))= K$, $$D(H_{\D})(\k G_C(D(H_{\D})))^+ = D(H)\{ g \gamma -1: \, g\gamma \in K\}.$$  
For a group $A$ and a subgroup $B \subset A$, let $$B^{\perp}=\{ f \in \widehat{A} : \, f(b)=1 \mbox{ for all } b \in B\}.$$ Note that $K^{\perp} \subset \widehat{G \times \Gamma} \simeq \Gamma \times G.$
\begin{prop}\label{quotientKR}
For $\beta \in G({H_{\D}}^*)=\Gamma$ and $g \in G$, the simple $D(H_{\D})$-module $H_{\D}\db g$ is also a $\overline{D(H_{D})}$-module via the quotient map, if and only if $\beta g \in K^{\perp}$.
\end{prop}
\begin{proof}
$H_{\D} \db g$ is a $\overline{D(H_{\D})}$-module, if and only if $f\gamma \cdot(h\db g) = h\db g$, for all $f\gamma \in K$ and $h \in H_{\D}$. Since $K \subset \mathcal{Z}(D(H_{\D}))$, if $f\gamma \in K$ then $ f\gamma \cdot (h \db g) = (f\gamma h) \cdot g = (h f \gamma) \cdot g = h \db ((f\gamma)\cdot g)$. Thus, $H_{\D} \db g$ is a $\overline{D(H_{\D})}$-module, if and only if $f\gamma \cdot g = g$, for all $f\gamma \in K$. Now $f\gamma \cdot g = f \db \gamma(g)g = \gamma(g) \beta(f)g$. And so, $H_{\D} \db g$ is a $\overline{D(H_{\D})}$-module, if and only if $  \gamma(g)\beta(f) =1$ for all $f \gamma \in K$; that is, if and only if, $ \beta g \in K^{\perp}$.
\end{proof}
\begin{lemma}
$K^{\perp} = \langle \mbox{\rm ev}_{\chi^{-1}a} \rangle $.
\end{lemma}
\begin{proof}
Since $K^{\perp} \simeq \widehat{\left(\frac{G \times \Gamma}{K}\right)}$, we have $|K^{\perp}|= | \frac{G \times \Gamma}{K} | = |\Im \mbox{\rm ev}_{\chi^{-1}a}|=|\mbox{\rm ev}_{\chi^{-1}a}|$; the last equality holding as $\Im \mbox{\rm ev}_{\chi^{-1}a}$ is cyclic (since it is a finite subgroup of $\k^{\times}$).  By the definitions of $K$ and $K^{\perp}$, $\mbox{\rm ev}_{\chi^{-1}a} \in K^{\perp}$, hence $K^{\perp} = \langle \mbox{\rm ev}_{\chi^{-1}a} \rangle $.
\end{proof}
It will be convenient to think of $K^{\perp}$ as a subgroup of $G \times \Gamma$ via the identification $\widehat{G \times \Gamma} \simeq  \widehat{G} \times \widehat{\Gamma} \simeq  \Gamma  \times G \simeq G \times \Gamma$. Under this identification, $K^{\perp} = \langle a\chi^{-1} \rangle$. 
\begin{rem}\label{simple-modules-quotient-KR} We can restate Proposition \ref{quotientKR} as follows:
the simple $D(H_{\D})$-modules that are also $\overline{D(H_{\D})}$-modules are of the form $H_{\D} \DOT_{(\chi^{-c})}a^c$, for $c=1, \dots,|a\chi^{-1}|.$
\end{rem}

Recall that $\mathcal{S}_{D(H_{\D})}$ denotes the set of isomorphism classes of simple $D(H_{\D})$-modules.
Combining Proposition \ref{tensor products}, Corollary \ref{1-dim-KR} and Proposition \ref{quotientKR}, we get that the map $$\Phi:\mathcal{S}_{\overline{D(H_{\D})}}  \times \mathcal{S}_{D(H_{\D})}^1 \to \mathcal{S}_{D(H_{\D})}$$ given by
$\Phi([U], [V]) =[ U\otimes V]$, is equivalent to the multiplication map $$\mu : K^{\perp} \times K \to G \times \Gamma, $$ under the identification of simple $D(H_{\D})$-modules with elements of $G \times \Gamma$.

\begin{thm}\label{decomposition-KR}
The map $\Phi$ as above is a bijection if and only if $\ell$ is odd and $\ell=M=N$.
\end{thm}

\begin{proof}
By the last remark, $\Phi$ is an bijection, if an only if $G\times \Gamma =  K^{\perp}\times K$; that is, if and only if $G \times \Gamma = K^{\perp}K$ and $K\cap K= \{1\}$. Now $|K^{\perp}| = | \frac{G \times \Gamma}{K} | = \frac{ |G \times \Gamma |}{|K|}$, and so $|K^{\perp}K| = \frac{ |K^{\perp} | |K|}{| K^{\perp} \cap K |}= \frac{ | G \times \Gamma|}{| K^{\perp} \cap K |}$. Hence $K^{\perp}K = G\times \Gamma$ if and only if $K^{\perp} \cap K = \{ 1 \}$. 

If $\ell=M=N$, then $|a|=|\chi|=\ell$ and so $ |a\chi^{-1}| = \ell$. Since $K^{\perp} \cap K \subset K^{\perp} = \langle a \chi^{-1} \rangle$, we have that $K^{\perp} \cap K = \langle (a\chi^{-1})^r \rangle$ for some $r \in \{ 1, \cdots ,  \ell \}$.
Since $(a\chi^{-1})^r \in K=\ker(\mbox{ev}_{\chi^{-1}a})$, $1= \mbox{\rm ev}_{\chi^{-1}a}\left( (a\chi^{-1})^r\right)= \left(\chi^{-1}(a)\right)^{2r}$ and so $\ell|\,2r$. If $\ell$ is odd, then $\ell|\,r$ and so $(a\chi^{-1})^r=1$, giving $K^{\perp} \cap K =\{1\}$.

Conversely, if $K^{\perp} \cap K =\{1\}$, let $n=|a\chi^{-1}|$. Then for all $r \in \{1, \cdots ,n -1 \}$, $(a\chi^{-1})^r \not\in K$. If either $M \neq \ell$ or $N \neq \ell$, then $ n>\ell$ and so $(a\chi^{-1})^{\ell} \not\in K$, which is a contradiction since $\mbox{\rm ev}_{\chi^{-1}a}((a\chi^{-1})^{\ell})=\chi^{-1}(a)^{2\ell} =1$. Hence, $\ell=M=N$.
If $\ell$ is even, then $(a\chi^{-1})^{\frac{\ell}{2}} \not\in K$, which is again a contradiction since $\mbox{\rm ev}_{\chi^{-1}a}((a\chi^{-1})^{\frac{\ell}{2}})= \chi^{-1}(a)^{\ell}=1$. Hence $\ell$ is odd.
\end{proof}
Next we describe the structure of $\overline{D(H_{\D})}$ under the hypothesis of the last theorem. 
\begin{prop}\label{pointed-u2}
If $\ell$ is odd and $\ell=N=M$, then $\overline{D(H_{\D})} \simeq \mathfrak{u}_{\theta}({\mathfrak{sl}}_2)$ as Hopf algebras, where $\theta = \chi (a)^{-\frac{1}{2}}$.
\end{prop}
\begin{proof}
Recall that $\mathfrak{u}_{\theta}({\mathfrak{sl}}_2) = \mathfrak{u}_{\theta,\theta^{-1}}({\mathfrak{sl}}_2)/ \langle (\w_1')^{-1} -\w_1 \rangle$. Since there is only one generator of each kind, we will omit the subindex 1. We have that $\mathfrak{u}_{\theta}({\mathfrak{sl}}_2)$ is generated by $e$, $f$ and $\w$, with relations:
$$ e^{\ell}=0=f^{\ell}, \; \w^{\ell} =1, \; \w e=\theta^2e\w, \; \w f= \theta^{-2}f\w \; \mbox{and} \; [e,f]=\frac{1}{\theta-\theta^{-1}}(\w - \w^{-1}).$$
In the proof of the previous proposition, we showed that if $\ell$ is odd and $\ell=N=M$, then $G\times \Gamma = \langle a \chi^{-1} \rangle K$, and so $\langle a\chi^{-1} \rangle$ is a complete set of representatives of the classes in $\frac{G\times \Gamma}{K}$.
Let $\psi : D(H_{\D}) \to \mathfrak{u}_{\theta}(\mathfrak{sl}_2)$ be the algebra map such that
\begin{itemize}
\item $\psi(g\gamma ) = \w^{-2c}$ if $ g\gamma \in (a\chi^{-1})^cK, \, \forall g\gamma \in G \times \Gamma$,
\item $\psi(\xi) = e$ and
\item $\psi(x) =(\theta - \theta^{-1})f$.
\end{itemize}
For $\psi$ to be defined, it must commute with the defining relations of $D(H_{D})$ (from Definition \ref{definition:HD} and Propositions \ref{definition:dualHD} and \ref{double-HD}). This is the case by the following calculations:
\begin{enumerate}

\item $\psi(x)\psi(g)= \chi(g) \psi(g)\psi(x) , \mbox{ for all } g \in G$:

Let $g \in G$; if $g \in (a\chi^{-1})^cK$, then $ga^{-c}\chi^c \in K$ and so $\chi^{-1}(g)\chi(a)^{2c}= 1$. Hence, $\chi(g)= \chi(a)^{2c}= \theta^{-4c}$ and
\begin{equation*}
\begin{split}
\psi(x)\psi(g) &=(\theta - \theta^{-1})f\w^{-2c}=(\theta - \theta^{-1})\theta^{-4c}\w^{-2c}f \\ 
&=\chi(g)\w^{-2c}(\theta - \theta^{-1})f =\chi(g)\psi(g)\psi(x).
\end{split}
\end{equation*}

\item $\psi(\xi)\psi(\gamma)= \gamma(a) \psi(\gamma)\psi(\xi) , \mbox{ for all } \gamma \in \Gamma$:

Let $\gamma \in \Gamma$, in a similar way as in the previous relation, it can be shown that if $\gamma \in (a\chi^{-1})^cK$, then $\gamma(a)= \chi(a)^{-2c}=\theta^{4c}$. We then have  
$$\psi(\xi)\psi(\gamma) = e \w^{-2c} = \theta^{4c} \w^{-2c} e= \gamma(a)\psi(\gamma)\psi(\xi).$$

\item $[\psi(x), \psi(\xi)] = \psi(a) -\psi(\chi)$:

To prove this, we first need to know the images of $a$ and $\chi$ under $\psi$. 
Since $\ell$ is odd, let $c \in \Z$ be such that $2c=1\mbox{ mod}\ell$. Then, $a=(a\chi^{-1})^c(a\chi)^c$, and since $a\chi \in K$, we have that \begin{equation*}\label{psi:a}\psi(a) = \w^{-2c}=\w^{-1}.\end{equation*}
Similarly, $\chi= (a\chi^{-1})^{-c}(a\chi)^{c}$ and so $\psi(\chi)= \w$.  
Now 
\begin{equation*}
\begin{split}
[\psi(x), \psi(\xi)]&= (\theta-{\theta}^{-1})[f,e]=-(\theta-{\theta}^{-1})[e,f] \\ & =-\frac{\theta-{\theta}^{-1}}{\theta-{\theta}^{-1}}(\w-\w^{-1})
=\w^{-1}-\w= \psi(a)- \psi(\chi).
\end{split}
\end{equation*}
\end{enumerate}
Clearly $\psi(x)^{\ell} = 0 = \psi (\xi)^{\ell}$ and $\psi(g)\psi(\gamma)= \psi(\gamma)\psi(g)$ for all $g\in G$ and $\gamma \in \Gamma$. The other relations follow in a similar way as 1 and 2 above.

Next we need to show that $\psi$ is a map of coalgebras. Group-like elements in $D(H_{\D})$ are mapped to group-like elements in $\mathfrak{u}_{\theta}(\mathfrak{sl}_2)$. Moreover,
\begin{eqnarray*}
\psi\otimes \psi (\com (x) ) &=& \psi \otimes \psi ( x \otimes a + 1 \otimes x)=(\theta - \theta^{-1})\left(f \otimes \w^{-1} + 1 \otimes f \right)\\
&=&(\theta - \theta^{-1})\com(f) = \com(\psi(x))
\end{eqnarray*}
and
\begin{eqnarray*}
\psi\otimes \psi (\com(\xi) ) = \psi \otimes \psi ( \chi \otimes \xi + \xi \otimes 1)=\left(\w \otimes e + e \otimes 1\right)=\com(e) = \com(\psi(\xi)).
\end{eqnarray*}
Therefore $\psi$ is a map of Hopf algebras.

Recall that $D(H_{\D})(\k K)^+=D(H_{\D})\left\{ k-1 \, : \, k \in K\right\}$. Note that $\psi(K)=\{1\}$ and so $\psi\left(\left\{ k-1 \, : \, k \in K\right\}\right) = \{0\}$. Therefore $D(H_{\D})(\k K)^+ \subset \ker(\psi)$ and the map $\psi$ induces a Hopf algebra map $\overline{\psi} : \overline{D(H_{\D})} \to \mathfrak{u}_{\theta}({\mathfrak{sl}}_2)$. Since $\ell$ is odd, $\langle \w \rangle = \langle \w^{-2} \rangle$, and so $\overline{\psi}$ is surjective.

By Remark \ref{schneider:quotient}, 
\begin{eqnarray*}
\dim(\overline{D(H_{\D})}) &=& \frac{\dim(D(H_{\D}))}{\dim(\k K)}= \frac{|G\times \Gamma|\ell^2}{|K|}=|K^{\perp}|\ell^2 = |\langle a\chi^{-1}\rangle| \ell^2 = \ell^3\\ &=&\dim(\mathfrak{u}_{\theta}({\mathfrak{sl}}_2)).\end{eqnarray*} Hence, $\overline{\psi}$ is an isomorphism.

\end{proof}

\begin{rem}\label{correspondence:modules}
Let $\mathfrak{b}'$ be (as in Chapter II) the subalgebra of $\mathfrak{u}_{\theta, \theta^{-1}}(\mathfrak{sl}_2)$ generated by $f$ and $\w'$ and $H= (\mathfrak{b'})^{\coop}$. Via the isomorphism $\overline{\psi}$ defined in the proof of Proposition \ref{pointed-u2}, a simple $D(H_{\D})$-module of the form $H_{\D}\DOT_{(\chi^{-c})}(a^c)$ is also a $\mathfrak{u}_{\theta}(\mathfrak{sl}_2)$-module. Explicitly, for $h \in \mathfrak{u}_{\theta}(\mathfrak{sl}_2)=\overline{\mathfrak{u}_{\theta,\theta^{-1}}(\mathfrak{sl}_2)}$ and $m \in H_{\D}\DOT_{(\chi^{-c})}(a^c)$, $h \cdot m = {\overline{\psi}}^{-1}(h) \cdot m$. Therefore, by Example \ref{u2}, as $\mathfrak{u}_{\theta}(\mathfrak{sl}_2)$-modules $H_{\D}\DOT_{(\chi^{-c})}(a^c) \simeq H \db (\w'^d)$ with $\beta(\w')=\theta^{-2d}$ for some $d \in \Z$. By analyzing the action of $\w'$ on both of these modules, it follows that $d=-c$.
Conversely, a simple $\mathfrak{u}_{\theta}(\mathfrak{sl}_2)$-module $H\db (w')^d$ becomes a simple $D(H_{\D})$-module via $\overline{\psi}$, and is isomorphic to $H_{\D}\DOT_{(\chi^{d})}(a^{-d})$ as $D(H_{\D})$-modules.
\end{rem} 

Next we study the reducibility of tensor products of simple $D(H_{\D})$-modules when $n=M=N$ is odd.

In \cite{radford} Radford used his construction to describe simple modules for the Drinfel'd Double of the Taft algebra, which is isomorphic to $\mathfrak{u}_{\theta, \theta^{-1}}(\mathfrak{sl}_2)$ when $\ell$ is odd ($\ell$ is the order of $\theta$).

Translating his result to our notation ($H=(\mathfrak{b}')^{\coop})$, generated by $\w'$ and $f$ and the corresponding relations) we have

\begin{prop}[Radford \cite{radford}]\label{radford-taft}
For $g=(\w')^c$ and $\beta: H \to \k$ an algebra morphism, let $r \geq 0$ be minimal such that $\beta(\w')=\theta^{2(c-r)}$. Then the simple $\mathfrak{u}_{\theta,\theta^{-1}}(\mathfrak{sl}_2)$-module $H\db g$ is $(r+1)-$dimensional with basis $\left\{ g, f \db g, \dots, f^r \db g \right\}$ and $f^{r+1} \db g =0.$
\end{prop}
In \cite{chen}, H-X. Chen studied the reducibility of tensor products of these simple modules; we translated Chen's result into Radford's notation:

\begin{prop}[Chen \cite{chen}]\label{chen}
Given $g=(\w')^c$, $g'=(\w')^{c'}$ in $G(H)$ and $\beta, \beta' \in G(H^*)$, let $r, r' \in \left\{ 0, \dots, \ell -1 \right\}$ be such that $\beta(\w')=\theta^{2(c-r)}$ and $\beta'(\w')=\theta^{2(c'-r')}$. Then the $\mathfrak{u}_{\theta,\theta^{-1}}(\mathfrak{sl}_2)$-module $H\db g \otimes H \DOT_{\beta'}g'$ is completely reducible if and only if $r+r' < \ell$. Moreover, let $$g_j= gg'(\w')^{-j}\quad \mbox{and} \quad \beta_j (\w')=\theta^{2j}\beta(\w')\beta'(\w');$$ if $r+r' < \ell$ then
$$H\db g \otimes H \DOT_{\beta'}g'\simeq \bigoplus_{j=0}^{\min(r,r')} H \DOT_{\beta_j} g_j.$$ 
If $r+r' \geq \ell$, let $t= r+r'-\ell +1$; then
$$\mbox{\rm Soc} \left( H\db g \otimes H \DOT_{\beta'}g'\right)\simeq \bigoplus_{j=\left[\frac{t+1}{2}\right]}^{\min(r,r')} H \DOT_{\beta_j} g_j.$$
\end{prop}

\begin{rem}\label{obs:u2}
Let $g= (\w')^c \in G(H)$; by Example \ref{u2}, if $H\db g$ is naturally a $\mathfrak{u}_{\theta}(\mathfrak{sl}_2)$-module, then $\beta=\beta_g$, {\it i.e.} $\beta(\w')=\theta^{-2c}=\theta^{2(c-2c)}$. Then the number $r$ from Proposition \ref{radford-taft} is $r=2c\mbox{ mod}\ell$, with $0\leq r < \ell$. We will denote such number by $r_c$. 
\end{rem}
We get the following corollary for simple $\mathfrak{u}_{\theta}(\mathfrak{sl}_2)$-modules:

\begin{cor}\label{chen:u2} Given $g=(\w')^c$ and $g'=(\w')^{c'}$ in $G(H)$. If $r_c + r_{c'} < \ell$ then 
$$H\DOT_{\beta_g}g \otimes H \DOT_{\beta_{g'}}g' \simeq \bigoplus_{j=0}^{\min(r_c,r_{c'})} H \DOT_{\beta_j} g_j,$$ as $\mathfrak{u}_{\theta}(\mathfrak{sl}_2)$-modules, where $g_j=gg'(\w')^{-j}$ and $\beta_j=\beta_{g_j}$.
\end{cor}

\begin{rem}
This last corollary is a particular case of a more general formula for simple modules for the non-restricted quantum group $U_{\theta}(\mathfrak{sl}_2)$, that appears as an exercise in \cite{bakalov-kirillov}.
\end{rem} 
 
We have an analogous result to Proposition \ref{chen} for $D(H_{\D})$-modules:

\begin{thm}\label{last-pointed} If $\ell=M=N$ is odd and $g\beta$, $g '\beta' \in G \times \Gamma=G(D(H_{\D}))$, let $c$ and $c' \in \Z$ be such that $(a^{-1}\chi)^c$ and $(a^{-1}\chi)^{c'}$ are representatives of the classes of  $g\beta$ and $g'\beta'$ in $G\times \Gamma/K$ respectively. Then the $D(H_{\D})$-module $ H_{\D}\db g \otimes H_{\D} \DOT_{\beta'}g' $ is completely reducible if and only if $r_{c} + r_{c'} < \ell$. Moreover, let $$ g_j=gg'a^{j}\quad \mbox{and} \quad \beta=\chi^{-j} \beta  \beta';$$if $r_{c} + r_{c'} < \ell$ then  
$$H_{\D}\db g \otimes H_{\D} \DOT_{\beta'}g' \simeq \bigoplus_{j=0}^{\min(r_{c},r_{c'})} H_{\D}\DOT_{\beta_j}g_j.$$
If $r_{c} + r_{c'} \geq \ell$, then $$\mbox{\rm Soc}\left(H_{\D}\db g \otimes H_{\D} \DOT_{\beta'}g'\right) \simeq \bigoplus_{j=\left[\frac{t+1}{2}\right]}^{\min(r_{c},r_{c'})} H_{\D}\DOT_{\beta_j}g_j,$$
where $t=r_{c} + r_{c'} - \ell +1$.
\end{thm}
\begin{proof}
Let $g_K\beta_K$ and $g'_K \beta'_K \in K$ such that $g\beta= (a^{-1}\chi)^{c}g_K\beta_K$ and $g'\beta'=(a^{-1}\chi)^{c'}g'_K \beta'_K $.
By Proposition \ref{decomposition-KR}, $H_{\D}\db g \simeq H_{\D} \DOT_{\chi^{c}} a^{-c} \otimes H_{\D}\DOT_{\beta_K}g_K$, the first factor in $\mathcal{S}_{\overline{D(H_{\D})}}$, and the second factor in $\mathcal{S}_{D(H_{D})}^1$. Similarly $H_{\D}\DOT_{\beta'} g' \simeq H_{\D} \DOT_{\chi^{c'}} a^{-c'} \otimes H_{\D}\DOT_{\beta'_K}g'_K$.
Then 
\begin{eqnarray*}
H_{\D}\db g \otimes H_{\D} \DOT_{\beta'}g' &\simeq& \left(H_{\D} \DOT_{\chi^{c}} a^{-c}\otimes H_{\D}\DOT_{\beta_K}g_K \right) \otimes \left(  H_{\D} \DOT_{\chi^{c'}} a^{-c'} \otimes H_{\D}\DOT_{\beta'_K}g'_K \right) \\
&\simeq&  \left(H_{\D} \DOT_{\chi^{c}} a^{-c} \otimes  H_{\D} \DOT_{\chi^{c'}} a^{-c'} \right) \otimes \left(H_{\D}\DOT_{\beta_K}g_K \otimes  H_{\D}\DOT_{\beta'_K}g'_K \right)\\
&\simeq&\left(H_{\D} \DOT_{\chi^{c}} a^{-c} \otimes  H_{\D} \DOT_{\chi^{c'}} a^{-c'} \right)  \otimes H_{\D} \DOT_{{}_{\beta_K  \beta'_K}}g_Kg'_K;
\end{eqnarray*}
the second isomorphism by symmetry of tensor products of modules for $D(H_{\D})$, and the third by Proposition \ref{tensor products}.
Let $\gamma, \gamma' :H \to \k$ be the algebra maps given by $\gamma(\w')= \theta^{-2c}$ and $\gamma'(\w')=\theta^{-2c'}$. If $r_c + r_{c'} < \ell$, we have the following isomorphisms of $\mathfrak{u}_{\theta}(\mathfrak{sl}_2)$-modules:
$$H_{\D} \DOT_{\chi^{c}} a^{-c} \otimes  H_{\D} \DOT_{\chi^{c'}} a^{-c'} \simeq
H \DOT_{\gamma}(\w')^{c}  \otimes  H \DOT_{\gamma'} (\w')^{c'} \simeq
\bigoplus_{j=0}^{\min(r,r')} H \DOT_{\beta_j} g_j,$$ where $g_j= (\w')^{c+c'-j}$and $\gamma_j (\w')=\theta^{-2(c+c'-j)}$, the first isomorphism following from Remark \ref{correspondence:modules} and the second from Corollary \ref{chen:u2}. Again by Remark \ref{correspondence:modules}, the $j^{\text th}$ summand of the last module is isomorphic to $H_{\D}\DOT_{\chi^{-c_j}}a^{c_j}$ as $D(H_{\D})$-modules, where $c_j= -(c+c'-j)$.
Then
$$H_{\D}\db g \otimes H_{\D} \DOT_{\beta'}g' 
\simeq
\left( \bigoplus_{j=0}^{\min(r,r')}H_{\D}\DOT_{\chi^{-c_j}}a^{c_j} \right) \otimes  H_{\D} \DOT_{{}_{\beta_K  \beta'_K}}g_Kg'_K 
\simeq
\bigoplus_{j=0}^{\min(r,r')}H_{\D}\DOT_{\gamma_j}g_j,$$ 
where $$g_j=a^{c_j}g_Kg'_K= a^{-c}g_Ka^{-c'}g'_Ka^j= gg'a^j$$ and
$$\gamma_j= \chi^{-c_j}\beta_K\beta'_K= \chi^c\beta_K \chi^{c'}\beta'_K \chi^{-j}=\beta \beta' \chi^{-j}.$$

If $r_c+r_{c'} \geq \ell$, we have
$$H_{\D}\db g \otimes H_{\D} \DOT_{\beta'}g' \simeq \left( H_{\D} \DOT_{\chi^{-c}} a^{c} \otimes  H_{\D} \DOT_{\chi^{-c'}} a^{c'} \right)_{\beta_K\beta'_{K}}$$ and so
$$\mbox{Soc}\left(H_{\D}\db g \otimes H_{\D} \DOT_{\beta'}g'\right) \simeq \left(\mbox{Soc}\left( H_{\D} \DOT_{\chi^{-c}} a^{c} \otimes  H_{\D} \DOT_{\chi^{-c'}} a^{c'} \right)\right)_{\beta_K\beta'_{K}}.$$ With a similar reasoning as before, we get that $$\mbox{Soc}\left( H_{\D} \DOT_{\chi^{-c}} a^{c} \otimes  H_{\D} \DOT_{\chi^{-c'}} a^{c'} \right) \simeq \bigoplus_{j=\left[\frac{t+1}{2}\right]}^{\min(r_{c}, r_{c'})}H \DOT_{\chi^{-c_j}}a^{c_j},$$ where $c_j=-(c+c'-j)$. Therefore
\begin{eqnarray*}
\mbox{Soc}\left(H_{\D}\db g \otimes H_{\D} \DOT_{\beta'}g'\right) &\simeq& \left(\bigoplus_{j=\left[\frac{t+1}{2}\right]}^{\min(r_{c}, r_{c'})}H \DOT_{\chi^{-c_j}} a^{c_j}\right)_{\beta_K\beta'_K}\\
& \simeq & \left( \bigoplus_{j= \left[\frac{t+1}{2}\right]}^{\min(r_{c}, r_{c'})}H \DOT_{\chi^{-c_j}}a^{c_j} \right) \otimes H \DOT_{{}_{\beta_K  \beta'_K}}g_Kg'_K \\
&\simeq& \bigoplus_{j=\left[\frac{t+1}{2}\right]}^{\min(r_{c},r_{c'})} H\DOT_{\beta_j}g_j,
\end{eqnarray*}
where $g_j= h_jg_Kg'_K=a^{-c-c'+j}g_Kg'_K= gg'a^j$ and $\beta_j=\gamma_j\beta_K  \beta'_K = \beta\beta'\chi^{-j}.$
\end{proof}

\begin{rem}
Following the notation of the proof of the last proposition, since $H_{\D}\db g \simeq H_{\D} \DOT_{\chi^{c}} a^{-c} \otimes H_{\D}\DOT_{\beta_K}g_K$ and $\dim ( H_{\D}\DOT_{\beta_K}g_K) =1$, we have $\dim (H_{\D} \db g) =\dim (H_{\D} \DOT_{\chi^{c}} a^{-c})$. Furthermore, $ H_{\D} \DOT_{\chi^{c}} a^{-c} \simeq H \DOT_{\gamma} (\w')^c $ as $\mathfrak{u}_{\theta}(\mathfrak{sl}_2)$-modules. Therefore $\dim (H_{\D} \db g) = \dim ( H \DOT_{\gamma} (\w')^c ) = r_c +1$ and similarly, $\dim (H_{\D} \DOT_{\beta'} g') = r_{c'} +1$. Then, the condition for complete reducibility of $H_{\D} \db \otimes H_{\D} \DOT_{\beta'} g'$ is $r_c + r_{c'} = \dim (H_{D} \db ) + \dim (H_{\D} \DOT_{\beta'} g')-2 < \ell$; that is $\dim (H_{D} \db ) + \dim (H_{\D} \DOT_{\beta'} g') \leq \ell +1$
\end{rem}

In \cite{e-g-s-t}, the authors studied the representation theory of the Drinfel'd double of a family of Hopf algebras that generalize the Taft algebra. In their case, the order of the generating group-like element need not be the same as the order of the root of unity. They give a similar decomposition of tensor products as in Theorem \ref{last-pointed}. Although the algebras $H_{\D}$ generalize their Hopf algebras, Theorem \ref{last-pointed} does not generalize their result since we require $|a|=|\chi(a)|$. However, since $G$ need not be cyclic, Theorem \ref{last-pointed} generalizes Chen's result for Taft algebras.

\end{document}